\documentclass[11pt,a4paper]{amsart}

\usepackage{graphicx}      

\usepackage{amsmath}
\usepackage{amssymb}
\usepackage{xcolor}
\usepackage{amsthm}
\usepackage{nicefrac}

\newtheorem{thm}{Theorem}
\newtheorem{cor}[thm]{Corollary}

\newtheorem{prop}[thm]{Proposition}
\newtheorem{assum}[thm]{Assumption}
\newtheorem{exmp}[thm]{Example}
\theoremstyle{definition}

\newtheorem{rem}[thm]{Remark}

\newcommand{\Lmu}{L^2(\mathbb{X})}

\newcommand{\LV}{\mathcal{L}_\mathbb{V}}
\newcommand{\tL}{\tilde{\mathcal{L}}_m}
\newcommand{\tC}{\tilde{C}_m}
\newcommand{\tA}{\tilde{A}_m}
\newcommand{\XX}{\operatorname{\mathbb{X}}}
\newcommand{\calL}{\mathcal L}  

\usepackage{tikz}
\usetikzlibrary{shapes.geometric, arrows}
\tikzstyle{box} = [rectangle, rounded corners, minimum width=4cm, minimum height=1cm,text centered, draw=black, fill=black!30]
\tikzstyle{midbox} = [rectangle, rounded corners, minimum width=2.62cm, minimum height=1cm,text centered, draw=black, fill=black!30]
\tikzstyle{arrow} = [thick,->,>=stealth]

\begin{document}
\title[Towards reliable data-based optimal and predictive control using eDMD]{Towards reliable data-based optimal and predictive control using extended DMD}
\author[M.\ Schaller, K.\ Worthmann, F.\ Philipp, S.\ Peitz, F.\ Nüske]{Manuel Schaller$^{1}$, Karl Worthmann$^{1}$ Friedrich Philipp$^{1}$, Sebastian Peitz$^{2}$ and Feliks Nüske$^{2,3}$}
	\thanks{}
	\thanks{$^{1}$Technische Universit\"at Ilmemau, Institute of Mathematics, Optimization-based Control group, Germany (e-mail: \{friedrich.philipp,manuel.schaller,karl.worthmann\}@tu-ilmenau.de).}
	\thanks{$^{2}$Paderborn University, Department of Computer Science, Data Science for Engineering, Germany, (e-mail: sebastian.peitz@upb.de).}
		\thanks{$^{3}$Max Planck Institute for Dynamics of Complex Technical Systems, Magdeburg, Germany, (e-mail: nueske@mpi-magdeburg.mpg.de)}

	\thanks{{\bf Acknowledgments: }F.\ Philipp was funded by the Carl Zeiss Foundation within the project \textit{DeepTurb---Deep Learning in and from Turbulence}. 
		K.\ Worthmann gratefully acknowledges support by the German Research Foundation (DFG; grants WO\,2056/6-1, WO\,2056/14-1). }.

\begin{abstract}
	While Koopman-based techniques like extended Dynamic Mode Decomposition are nowadays ubiquitous in the data-driven approximation of dynamical systems, quantitative error estimates were only recently established. To this end, both sources of error resulting from a finite dictionary and only finitely-many data points in the generation of the surrogate model have to be taken into account. We generalize the rigorous analysis of the approximation error to the control setting while simultaneously reducing the impact of the curse of dimensionality by using a recently proposed bilinear approach. In particular, we establish uniform bounds on the approximation error of state-dependent quantities like constraints or a performance index enabling data-based optimal and predictive control with guarantees.
	
\end{abstract}

\maketitle

\section{Introduction}
	While optimal and predictive control based on models derived from first principles is nowadays well established, data-driven control design is becoming more and more popular. %
	We present an approach via extended Dynamic Mode Decomposition (eDMD) using the Koopman framework to construct a data-driven surrogate model suitable for optimal and predictive control. 
	
\begin{figure}[htb]
		\centering
			{\begin{tikzpicture}[node distance=5.2cm]
				\node (init) {};
				\node (start) [midbox, above right of=init, node distance=0.55cm] {observable $\varphi$}; 
				\node (mid) [midbox, right of=start] {$\mathcal{K}^t \varphi$};
				\node (rig) [midbox, below of=mid,node distance=2cm] {$(\mathcal{K}^t \varphi)(x_0)$};
				\draw [arrow] (start) -- node[anchor=south] {\small Koopman} (mid);
				\draw [arrow] (mid) -- node[anchor=west] {\small evaluation} (rig);
				\node (start2) [midbox, below left of=init, node distance=0.6cm] {initial state $x_0$}; 
				\node (mid2) [midbox, below of=start2, node distance=2cm] {$x(t;x_0)$};
				\node (rig2) [midbox, right of=mid2] {$\varphi(x(t;x_0))$};
				\draw [arrow] (start2) -- node[anchor=east] {\small ODE} (mid2);
				\draw [arrow] (mid2) -- node[anchor=south] {\small evaluation} (rig2);
		\end{tikzpicture}}
		
		\caption{Schematic sketch of the Koopman framework: Instead of first propagating the ODE and then evaluating the observable, the observable is propagated and then evaluated at the initial state.}
		\label{fig:sketch}
	\end{figure}
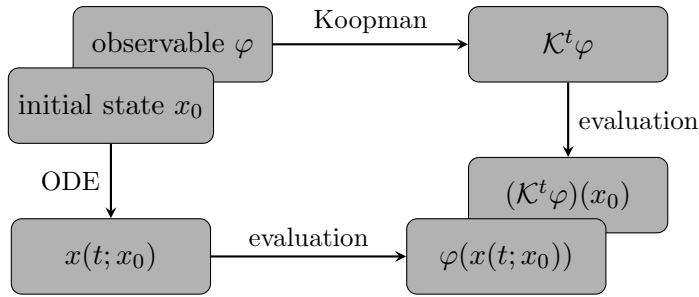

	The Koopman framework provides the theoretical foundation for data-driven approximation techniques like eDMD, see \cite[Chapters~1 and~8]{mauroy2020koopman}: %
	Using the Koopman 
	semigroup~$(\mathcal{K}^t)_{t \geq 0}$ or, equivalently, the Koopman generator~$\mathcal{L}$, 
	observables~$\varphi$ (real-valued $L^2$-functions of the state) 
	can be propagated forward-in-time via
	\begin{equation}\nonumber
		\mathcal{K}^t\varphi = \mathcal{K}^0 \varphi + \mathcal{L}\int_0^t  \mathcal{K}^s\varphi \,\mathrm{d}s.
	\end{equation}
	Propagating the observable along the linear Koopman operator via $\mathcal{K}^t\varphi$ and evaluating the result at a state~$x_0$ provides an alternative to calculating the solution~$x(t;x_0)$ of the underlying Ordinary Differential Equation (ODE) and then evaluating the observable as depicted in Figure~\ref{fig:sketch}.
	
	In the analysis of the eDMD-based approximation $(\tilde{\mathcal{K}}^t)_{t\geq 0}$ of the  Koopman semigroup $(\mathcal{K}^t)_{t\geq 0}$, two sources of error have to be taken into account: The projection and the estimation 
	error. %
	First, a dictionary is chosen, which consists of finitely-many observables $\psi_1,\ldots,\psi_N$ and, thus, spans a finite-dimensional subspace~$\mathbb{V}$. Since the eDMD-based surrogate model is constructed on~$\mathbb{V}$, a \textit{projection error} occurs. Second, only a finite number of data points $x_1,\ldots,x_m$ is used to generate the surrogate model, which induces an additional \textit{estimation error} on~$\mathbb{V}$. 
	Whereas the convergence of the eDMD-based approximation to the Koopman semigroup in the infinite-data limit, i.e., for $N$ and $m$ tending to infinity, was shown in~\cite{KM18b}, error bounds 
	for a finite dictionary and finite data 
	depending on~$N$ and~$m$ were derived 
	in~\cite{ZhanZuaz21} and~\cite{NuskPeit21} for identically-and-independently distributed (i.i.d.) data 
	for the estimation 
	step. While also the projection error is analyzed in 
	the former reference, the latter covers the estimation error even for stochastic differential equations and	ergodic sampling.

	We consider 
	the nonlinear control-affine ODE 
	\begin{equation}\label{eq:ode}
		\dot{x}(t) = g_0(x(t)) + \sum\nolimits_{i=1}^{n_c} g_i(x(t))\, u_i(t)
	\end{equation}
 with locally Lipschitz-con\-ti\-nu\-ous vector fields~$g_0, g_1, \ldots, g_{n_c}: \mathbb{R}^n \rightarrow \mathbb{R}^n$ and subject to the initial condition $x(0)=x_0$.
	Further, we impose the control constraints $u(t) \in \mathbb{U}$ for some compact, convex, and nonempty set $\mathbb{U} \subset \mathbb{R}^{n_c}$ and define, for $T > 0$, the set of \textit{admissible control functions} by
	\begin{equation}\label{eq:admissible_controls}
		\mathcal{U}_T(x_0) \triangleq 
		\left\{ u: [0,T] \rightarrow \mathbb{R}^{n_c} \left| %
		\begin{array}{l}
			\text{$u$ measurable} \\
			\exists!\,x(\cdot;x_0,u) \\
			u(t) \in \mathbb{U}, t \in [0,T]
		\end{array}
		\right. \right\},
	\end{equation}
	where $x(t;x_0,u)$ denotes the unique solution 
	at time~$t \geq 0$.
	
	\cite{Proctor2016} as well as~\cite{KM18a} proposed a method to predict control systems within the Koopman framework. To this end, the state is 
	augmented by the control variable. Then, a linear surrogate model depending on the extended state is generated using eDMD.
	Other popular methods are given by, e.g., using a coordinate transformation into Koopman eigenfunctions \cite{KKB21} or a component-wise Taylor series expansion~\cite{MCTM21}.
	In this work, however, we 
	use the bilinear approach, exploiting the control-affine structure of~\eqref{eq:ode} as suggested, e.g., in \cite{WillHema2016,Sura2016,Peitz2020}, 
	for which estimation error estimates were derived in~\cite{NuskPeit21}. The advantages of this approach are twofold. 
	First, one can observe a superior performance when considering nonlinear 
	systems with a control-state coupling, 
	which we briefly showcase in Example~\ref{ex:duffing}. Second, as the state dimension is not augmented, the data-requirements are less demanding. In particular, 
	the curse of dimensionality is alleviated in the multi-input case in comparison to the previously proposed state-augmentation. 

\noindent The probabilistic bounds 
	on the estimation error for the propagated observable derived in~\cite{NuskPeit21} depend on the 
	control function. However, for 
	optimal and 
	predictive control, it is 
	essential to derive uniform 
	estimates. Hence, our first key 
	contribution is to establish a bound 
	in Section~\ref{sec:Uniform}, which uniformly holds for all control functions on the prediction horizon.
	Our 
	second key 
	contribution 
	is the additional estimation of the projection error using a dictionary consisting of only finitely-many observables using 
	techniques well-known for finite-element methods in Section~\ref{Sec:projection}, see~\cite{Brae1997,QuarVall2008}. The derived bound decays with increasing size of the dictionary. In conclusion and to the best of the authors' knowledge, this is the first rigorous \textit{finite-data} error estimate for the eDMD-based prediction for nonlinear control systems taking into account both sources of errors, i.e., the projection \textit{and} the approximation error.

	The paper is organized as follows: In Section~\ref{Sec:Recap}, we briefly recap eDMD and the bilinear surrogate model obtained for control-affine control systems. Section~\ref{sec:Uniform} is devoted to rigorous error bounds on the estimation 
	error---uniform w.r.t.\ the control, while the projection error is considered in Section~\ref{Sec:projection}. Then, the application of the derived bounds in optimal and predictive control is discussed in Section~\ref{sec:control} before conclusions are drawn in Section~\ref{sec:conclusions}.

\section{Koopman generator and Extended DMD}\label{Sec:Recap}

In this section, we recap the extended Dynamic Mode Decomposition (eDMD) 
as an established methodology to generate 
a data-based surrogate model 
for 
the Koopman operator 
or its generator 
to approximately describe the dynamics of observables along the flow of 
the control-affine 
system~\eqref{eq:ode}, 
see~\cite{BBKK21,Mez05}. 

\subsection{eDMD for autonomous systems}\label{subsec:edmd}

In this subsection, 
we introduce the data-based finite-dimensional approximation of the Koopman generator and the corresponding Koopman operator for autonomous systems using eDMD, i.e., setting $u(t) \equiv \bar{u}\in \mathbb{U}$, see, e.g., \cite{WILLIAMS2015} and defining 
$\dot{x}(t) = f(x(t))$ by $f(x) = g_0(x) + \sum_{i=1}^{n_c} g_i(x) \bar{u}_i$. We consider this dynamical system on a compact 
set $\mathbb{X} \subsetneq \mathbb{R}^n$. 
For initial value $x_0 \in \mathbb{X}$, the Koopman 
semigroup acting on square-integrable 
measurable functions $\varphi \in L^2(\XX)$ is defined by $(\mathcal{K}^t \varphi)(x_0) = \varphi(x(t;x_0))$ on the maximal 
interval of existence of 
$x(\cdot;x_0)$. The corresponding Koopman generator $\mathcal{L}:D(\mathcal{L}) \subset L^2(\XX) \to L^2(\XX)$ is defined as 
\begin{align}\label{e:def_generator}
    \mathcal{L}\varphi := \lim_{t\to 0} \frac{(\mathcal{K}^t - \operatorname{Id})\varphi}{t}.
\end{align}
Hence, $z(t) = \mathcal{K}^t\varphi\in L^2(\XX)$ solves 
the Cauchy problem $\dot{z}(t) = \mathcal{L}z(t)$, $z(0) = \varphi \in D(\mathcal{L})$. 

For a 
dictionary of 
observables $\psi_1,\ldots,\psi_N\in D(\calL)$, we consider the finite-dimensional subspace 
\begin{align*}
    \mathbb{V} := \operatorname{span}\{ \psi_j, j = 1,\ldots,N \} \subset D(\mathcal{L}).
\end{align*} 
The orthogonal projection onto~$\mathbb{V}$ and the Galerkin projection of the Koopman generator are denoted by~$P_\mathbb{V}$ and $\LV: ={P}_{\mathbb{V}} \mathcal{L}|_{\mathbb{V}}$, resp. %
Along the lines of \cite{Klus2020}, we have the representation $\LV = C^{-1}A$ with $C, A \in \mathbb{R}^{N\times N}$,
\begin{align*}
    C_{i,j}=\langle\psi_i,\psi_j\rangle_{\Lmu} \quad\text{and}\quad A_{i,j} =\langle\psi_i,\mathcal{L}\psi_j\rangle_{\Lmu}.
\end{align*}
For data points $x_1,\ldots,x_{m} \in \mathbb{X}$ and the matrices
\begin{align*}
    \Psi(X) &:= \left( \left.   \left(\begin{smallmatrix}
        \psi_1(x_1)\\
        :\\
        \psi_N(x_1)
    \end{smallmatrix}\right)\right| \ldots \left| \left(\begin{smallmatrix}
        \psi_1(x_{m})\\
        :\\
        \psi_N(x_{m})
    \end{smallmatrix}\right)\right. \right)\\
    \mathcal{L}\Psi(X) &:= \left( \left. \left(\begin{smallmatrix}
        (\mathcal{L}\psi_1)(x_1)\\
        :\\
        (\mathcal{L}\psi_N)(x_1)
    \end{smallmatrix}\right)\right| \ldots \left| \left(\begin{smallmatrix}
        (\mathcal{L}\psi_1)(x_{m})\\
        :\\
        (\mathcal{L}\psi_N)(x_{m})
    \end{smallmatrix}\right)\right. \right),
\end{align*}
$(\mathcal{L}\psi_j)(x_i) = \langle f(x_i), \nabla \psi_j(x_i) \rangle$, define $\tC, \tA \in \mathbb{R}^{N\times N}$ by 
\begin{align*}
    \tC = \tfrac{1}{m} \Psi(X) \Psi(X)^\top \quad\text{and}\quad \tA = \tfrac{1}{m} \Psi(X) \mathcal{L}\Psi(X)^\top
\end{align*}
to obtain the empirical, i.e., purely data-based, estimator $\tL = \tC^{-1} \tA$ for the Galerkin projection $\LV$.

\subsection{Bilinear surrogate control system}\label{subsec:bilin}

We briefly sketch the main steps of the bilinear surrogate modeling approach as presented in \cite{WillHema2016,Sura2016, Peitz2020}, for which a finite-data bound on the estimation error was given in \cite{NuskPeit21}.
Since control affinity of the system is inherited by the Koopman generator, for $u \in L^{\infty}([0,T],\mathbb{R}^{n_c})$, we set 
\begin{align}
\label{e:bilin1}
    \calL^u(t) = \calL^0 + \sum_{i=1}^{n_c}u_i(t)\left(\calL^{e_i} - \calL^0\right),
\end{align}
where $\mathcal{L}^{e_i}$, $i \in \{0,\ldots,n_c\}$, is the Koopman generator for the autonomous system with constant control $\bar{u} = e_i$, where $e_0 = 0$. Then, we can describe the time evolution of an observable function $\varphi\in L^2(\mathbb{X})$ via the bilinear system \begin{align}
\label{eq:bilin_real}
    \dot{z}(t) = \calL^u(t)z(t),\qquad z(0)=\varphi,
\end{align}
where we omitted the control argument in $z(t)=z(t;u)$ for the sake of brevity. The propagated observable can then
be evaluated for an initial state~$x_0$ via $z(t;u)(x_0)$, cp.~Figure~\ref{fig:sketch}. 
The projection of~\eqref{e:bilin1} onto~$\mathbb{V}$, spanned by a finite dictionary, is given by $\LV^u (t) := \LV^0 + \sum_{i=1}^{n_c}u_i(t)\big(\LV^{e_i}-\LV^0\big)$; analogously to Subsection~\ref{subsec:edmd}. Hence, the propagation of an observable $\varphi \in L^2(\mathbb{X})$ projected onto~$\mathbb{V}$ is given by
\begin{align}\label{eq:bilin_proj}
    \dot{z}_\mathbb{V}(t) = \LV^u (t)z_\mathbb{V}(t), \qquad z(0) = \mathbb{P}_\mathbb{V}\varphi.
\end{align}
The corresponding approximation by means of eDMD using $m$~data points is defined analogously via
\begin{align}\label{eq:nonauto_approx}
\tL^u (t) &:= \tL^0 + \sum_{i=1}^{n_c}u_i(t)\big(\tL^{e_i}-\tL^0\big),
\end{align}
where 
$\tL^{e_i}$ are eDMD-based approximations of $\mathcal{L}^{e_i}_\mathbb{V}$. 
Then, the corresponding 
data-based surrogate model 
reads
\begin{align}
\label{eq:bilin_surr}
        \dot{\tilde{z}}_{m}(t) = \tL^u (t)\tilde{z}_{m}(t), \qquad \tilde{z}_{m}(0) = \mathbb{P}_\mathbb{V}\varphi.
\end{align}

Let us highlight that, contrary to the popular DMD with control (DMDc) approach \cite{Proctor2016,KM18a}, which yields linear surrogate models of the form $Ax+Bu$, 
numerical simulation studies indicate that bilinear surrogate models are better suited if control and state are \textit{coupled}, see Example~\ref{ex:duffing}.  
Another key feature of the bilinear approach is that the state-space dimension is not augmented by the number of 
inputs, which alleviates 
the curse of dimensionality in comparison to DMDc.
\begin{exmp}\label{ex:duffing}
    We briefly present an example with a Duffing oscillator, cf.\ \cite[Section 4.2.1]{NuskPeit21} for more details, using the bilinear 
    approach to showcase its superior performance compared to DMDc if state and control are coupled. To this end, consider the dynamics
    \begin{equation}\label{eq:Duffing_NonlinCon}
        \dot{x} 
        = \begin{pmatrix}
            x_2 \\ -\delta x_2 - \alpha x_1 - 2\beta x_1^3 u
        \end{pmatrix}, \quad x(0) = x_0,
    \end{equation}
    with $\alpha = -1$, $\beta = 1$, $\delta = 0$. 
    Figure~\ref{fig:Duffing_Prediction} shows the prediction accuracy for $m=100$ and the dictionary $\{\psi_j\}_{j=1}^{N}$ consisting of monomials with maximal degree five. 
    We observe an excellent agreement for the bilinear surrogate model for more than one second, whereas eDMDc yields a large error of approximately $10\%$ from the start and becomes unstable almost immediately. 
\begin{figure}[!ht]
    \centering
    \includegraphics[width=.8\columnwidth]{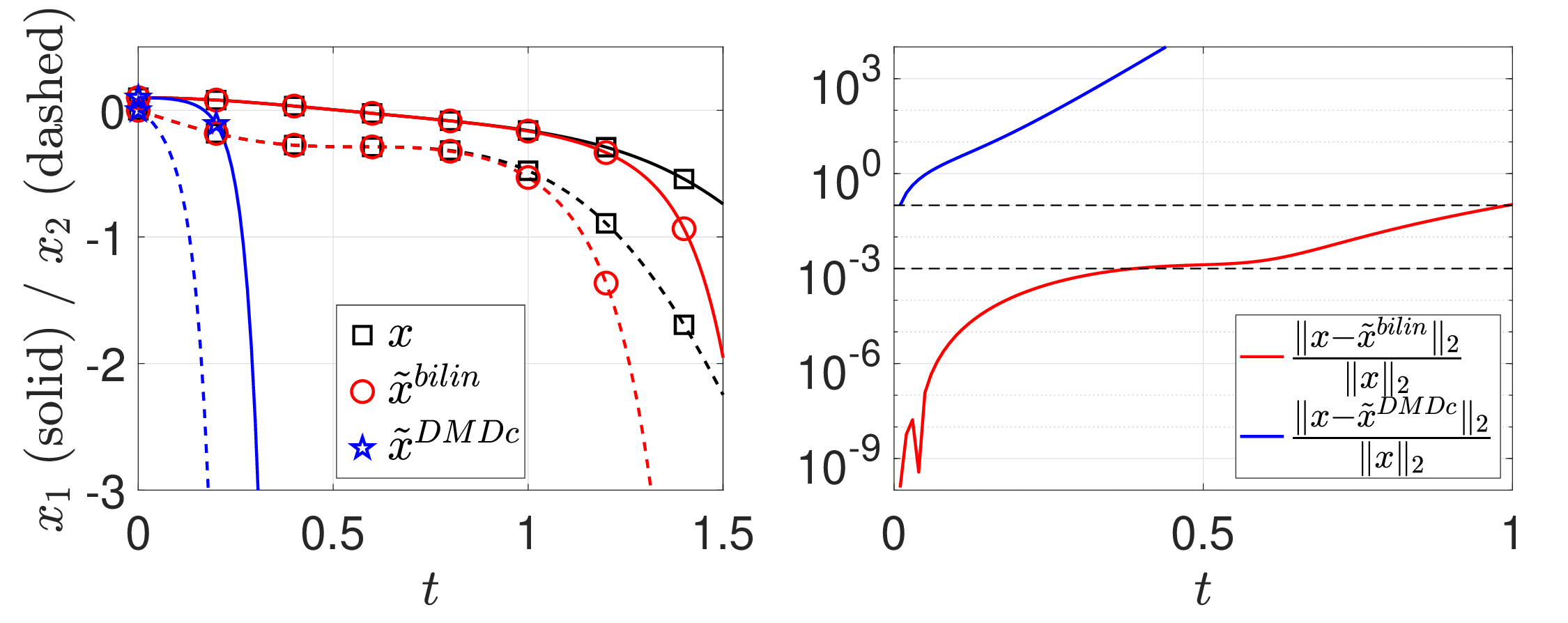}
    \caption{Comparison of the bilinear and the DMDc approach for~\eqref{eq:Duffing_NonlinCon} and a sinusoidal control input.}
\label{fig:Duffing_Prediction}
\end{figure}
\end{exmp}

\section{Estimation error: Uniform 
bounds} 
\label{sec:Uniform}

In this section, we derive an error bound that is uniform in the control~$u$ with values in the compact set $\mathbb{U}$ and, thus, refine the error bound of~\cite{NuskPeit21}. To this end, we require the following standard assumption.
\begin{assum}\label{as:main}
    Assume that the data, for each autonomous system with constant control $u \equiv e_i$, $i \in \{ 0,\ldots,n_c\}$, is sampled i.i.d.\ on $\mathbb{X}$ w.r.t.\ the Lebesgue measure. 
\end{assum}

We combine error bounds on the 
autonomous systems corresponding to $u \equiv e_i$, $i \in \{0,\ldots,n_c$\}, and exploit the control-affine structure of~\eqref{eq:ode} to derive the following error bound, which is an extension of our previous work by incorporating control constraints and providing a uniform bound independently of the chosen control function. 
\begin{thm}\label{cor:composite}
    Suppose that Assumption~\ref{as:main} holds and $\mathbb{U} \subset \mathbb{R}^{n_c}$ is bounded. 
    Then, for error bound $\varepsilon > 0$ and 
    probabilistic tolerance $\delta \in (0,1)$, the probabilistic error bound
    \begin{align}\label{eq:estimation_error}
       \mathbb{P}\big(\|\LV^u(t) - \tL^u(t)\|_F \leq \varepsilon\big) \geq 1-\delta \qquad\forall\,t \geq 0
    \end{align}	
    holds for all measurable control functions $u: [0,\infty) \rightarrow \mathbb{U}$ 
    if $m \geq \underbar{m}=\mathcal{O} ( \nicefrac {N^2} {\varepsilon^2 \delta} )$ holds for the number of data points, where $\| \cdot \|_F$ denotes the Frobenius norm. 

\end{thm}
\begin{proof}
   Invoking boundedness of $\mathbb{U}$, we set $\tilde{\delta}:=\nicefrac {\delta}{3(n_c+1)}$ and $\bar{\varepsilon}:=\nicefrac {\varepsilon}{(n_c+1)\left(1  + \max_{u \in \mathbb{U}} \sum_{i=1}^{n_c}|u_i| 
        \right)}$.
    For $k \in \{0,\ldots,n_c\}$, let the matrix $A^{(k)} \in \mathbb{R}^{N \times N}$ be defined by $\left(A^{(k)}\right)_{i,j} = \langle \psi_i,\mathcal{L}^{e_k} \psi_j \rangle_{L^2(\mathbb{X})}$ and set
    \begin{align*}
        \tilde{\varepsilon}_k &= \min\left\{1,\frac{1}{\|A^{(k)}\|\|C^{-1}\|}\right\} \cdot \frac {\|A^{(k)}\| \bar{\varepsilon}}{2\|A^{(k)}\|\|C^{-1}\| + \bar{\varepsilon}}.
    \end{align*}
    Then, choose a number of data points $\underbar{m}\in \mathbb{N}$ such that
        \begin{align}\label{eq:mindata}
            \underbar{m} \geq \max_{k=0,\ldots,n_c}\frac{N^2}{\tilde{\delta}\tilde{\varepsilon}_k^2} \max\left\{\|\Sigma_{A^{(k)}} \|^2_F, \|\Sigma_{C}\|^2_F\right\}
    \end{align}
    where $\Sigma_{A^{(k)}}$ and $\Sigma_{C}$ are variance matrices defined via 
    \begin{align*}
        \left(\Sigma_{A^{(k)}}\right)_{i,j}^2 &=  \int_\mathbb{X} \psi_i^2(x) \langle g_0(x) + g_k(x), \nabla \psi_j(x)\rangle^2\,\mathrm{d}x
        \\&\qquad \qquad - \left(\int_\mathbb{X} \psi_i(x) \langle g_0(x) + g_k(x), \nabla \psi_j(x) \rangle \,\mathrm{d}x\right)^2,\\
        \left(\Sigma_{C}\right)_{i,j}^2 & = \int_\mathbb{X} \psi_i^2(x) \psi_j^2(x)\,\mathrm{d}x - \left(\int_\mathbb{X} \psi_i(x)\psi_j(x)\,\mathrm{d}x\right)^2.
    \end{align*}
    Using $m \geq \underbar{m}$ data points, cp.~\eqref{eq:mindata}, we obtain probabilistic error estimates for the 
    generators $\tilde{\mathcal{L}}_m^{e_i}$, $i \in \{0,\ldots,n_c\}$, via \cite[Theorem 12]{NuskPeit21}:
    \begin{align}\label{eq:estimation_error_k}
        \mathbb{P}\left(\|\mathcal{L}_\mathbb{V}^{e_i}-\tilde{\mathcal{L}}^{e_i}_m\| \leq \bar{\varepsilon}\right) \geq 1-\tfrac{\delta}{n_c+1}.
    \end{align}    
    Rewriting $\mathcal{L}_\mathbb{V}^u (t) -\tilde{\mathcal{L}}_m^u(t)$ as
    \begin{align*}
        \Big( 1-\sum_{i=1}^{n_c}u_i(t) \Big) \left(\calL_\mathbb{V}^0- \tL^0\right) + \sum_{i=1}^{n_c}u_i(t) \left(\calL_\mathbb{V}^{e_i}-\tL^{e_i}\right),
    \end{align*}
    the desired error bound~\eqref{eq:estimation_error} can straightforwardly be derived based on the error bounds for the individual generators~\eqref{eq:estimation_error_k} analogously to \cite[Proof of Theorem 17]{NuskPeit21}. 
\end{proof}
Having a probabilistic bound for the estimation error on the projected non-autonomous generator at hand, a bound on the resulting trajectories of observables can be derived using Gronwall's inequality.
\begin{cor}
\label{cor:trajectories}
    Suppose that Assumption~\ref{as:main} holds and $\mathbb{U} \subset \mathbb{R}^{n_c}$ is bounded. 
    Let $T,\varepsilon > 0$, $\delta\in (0,1)$, and $z^0\in \mathbb{V}$ be given. Then, there is a number $\underbar{m} = \mathcal{O}( \nicefrac{N}{\varepsilon^2\delta} )$ of data points such that, for any $m \geq \underbar{m}$, the solutions $z,\tilde{z}_m$ of
    \begin{align*}
        &&\dot{z}(t) &= \LV^u(t)z(t), &&z(0)=z^0\\
        &&\dot{\tilde{z}}_m(t) &= \tL^u(t)\tilde{z}_m(t), &&\tilde{z}_m(0)={z}^0
    \end{align*}
satisfy
\begin{align*}
    \mathbb{P}\big( \|z(t)(x_0)-\tilde{z}_m(t)(x_0)\| \leq \varepsilon\big) \geq 1-\delta
\end{align*}
for all $x_0\in \mathbb{X}$, measurable control functions $u:[0,T] \rightarrow \mathbb{U}$ 
and $t\geq 0$ such that $x(s;x_0,u)\in \XX$ holds for all $s\in [0,t]$.
\end{cor}
\begin{proof}
The proof follows by straightforward modifications of \cite[Proof of Corollary 18]{NuskPeit21} using the uniform data requirements of Proposition~\ref{cor:composite}.
\end{proof}

Note that our approach to approximate the generator only requires the state to be contained in $\mathbb{X}$ up to any arbitrary small time $t>0$ to be able to define the generator as in \eqref{e:def_generator}. Then, in order to obtain error estimates for arbitrary long time horizons when going to a control setting, we have to ensure that the state trajectories remain in the set~$\mathbb{X}$ by means of our chosen control function. %
Besides a controlled forward-invariance of the set~$\mathbb{X}$, this can be ensured by choosing an initial condition contained in a suitable sub-level set of the optimal value function of a respective optimal control problem, see, e.g., \cite{BoccGrun14} or \cite{EsteWort20} for an illustrative application of such a technique in showing recursive stability of Model Predictive Control (MPC) without stabilizing terminal constraints for discrete- and continuous-time systems, respectively.

\section{Finite-data error bound for the approximation error}\label{Sec:projection}

In this section, we present our main result---a novel finite-data error bound for the full approximation error taking both estimation and projection error into account and, thus, generalizing~\cite[Proposition 5.1]{ZhanZuaz21} to non-autonomous and control systems. 

If the dictionary $\mathbb{V}$ 
forms a Koopman-invariant subspace, Corollary~\ref{cor:trajectories} directly yields an estimate for the observables, 
as the original system~\eqref{eq:bilin_real} and the projected system~\eqref{eq:bilin_proj} coincide.
If this is not the case, one further has to analyze the error resulting from projection onto the dictionary~$\mathbb{V}$. 
To this end, we choose a dictionary of finite elements.
\begin{assum}\label{ass:finite_elements}
    Suppose that the set $\mathbb{X}$ is compact and has a Lipschitz boundary~$\partial \mathbb{X}$. Further, let a regular, uniform triangulation of~$\mathbb{X}$ with meshsize $\Delta x>0$ be given. Further, let $\psi_i$ denote the (usual) linear hat function for the node~$x_i$, $i \in \{1,\ldots,N\}$, such that $\psi_i(x_j) = \delta_{ij}$ holds, where the latter is the Kronecker symbol.  
\end{assum}
The meshsize~$\Delta x$ might, e.g., be the incircle diameter of each cell. We point out that the size of the dictionary of finite elements is proportional to $\tfrac{1}{\Delta x^d}$ and refer to~\cite{QuarVall2008} and the references therein for details on finite elements. Furthermore, we emphasize that the dictionary~$\mathbb{V}$ consisting of the finite-elements functions may be further enriched by additional observables representing quantities of particular interest like state constraints or state-dependent stage costs. 

\begin{thm}\label{thm:fem}
    Suppose that Assumption~\ref{as:main} holds and that $\mathbb{U} \subset \mathbb{R}^{n_c}$ is bounded. 
    Let an observable $\varphi\in \mathcal{C}^2(\XX,\mathbb{R})$, an error bound $\varepsilon > 0$, a probabilistic tolerance $\delta \in (0,1)$, and a time horizon~$T > 0$ be given. 
    Then, if the dictionary consists of finite elements in accordance with Assumption~\ref{ass:finite_elements}, there is a mesh size $\Delta x=\mathcal{O}(\varepsilon)$ and a required 
    amount of data 
    $\underbar{m} = \mathcal{O} ( \nicefrac {1} {\varepsilon^{2+2d}\delta})$
    such that, for $\tilde{z}_{m}(0)=P_{\mathbb{V}}\varphi$, the probabilistic error bound
    \begin{align*}
        \mathbb{P} \left(\|\varphi(x(t;\cdot,u)) - \tilde{z}_{m}(t;\cdot,u)\|_{L^2(\mathcal{A}(t))} \leq \varepsilon \right) \geq 1-\delta
    \end{align*}
    holds for all measurable control functions $u:[0,T] \rightarrow \mathbb{U}$ and the data-based prediction using the bilinear surrogate dynamics~\eqref{eq:bilin_surr} generated with $m \geq \underbar{m}$ data points, where $\mathcal{A}(t) = \left\{ x_0 \in \XX\,|\, x(s;x_0,u) \in \XX\,\forall s \in [0,t] \right\}$.
\end{thm}
\begin{proof}
    First, we have $\varphi(x(t;x_0,u)) = z(t;u)(x_0)$, where $z$ solves~\eqref{eq:bilin_real}, i.e., 
    using $\mathcal{L}^0\varphi = g_0 \cdot \nabla \varphi$, $\mathcal{L}^{e_i}\varphi = (g_0 + g_i)\cdot \nabla \varphi$ and abbreviating $z(t)=z(t;u)$,
    \begin{align*}
        \dot{z}(t) = \mathcal{L}^{u(t)} z(t) &=\Big( \calL^0 + \sum_{i=0}^{n_c}u_i(t)\left(\calL^{e_i} - \calL^0\right) \Big) z(t)= \Big(g_0 + \sum_{i=0}^{n_c}u_i(t) g_i 
        \Big) \cdot \nabla z(t).
    \end{align*}
    This can be viewed as a linear transport equation 
    \begin{align} \label{eq:transport_z}
        \tfrac{\mathrm{d}}{\mathrm{d}t} z(t) = a(t,\cdot)\cdot\nabla z(t), \quad z(0) = \varphi, 
    \end{align}
    along the time- and space-dependent vector field 
    $$a(t,x) := g_0(x) + \sum_{i=0}^{n_c}u_i(t)\left(g_i(x) - g_0(x)\right).
    $$
    Since $\mathbb{X}$ is compact, $\mathbb{U}$ is bounded, and $g_i$, $i \in \{0,\ldots,n_c\}$, are continuous, there are $\underline{a},\overline{a}\in \mathbb{R}$ such that $\underline{a} \leq a(t,x) \leq \overline{a}$ for 
    a.e.\ $t \in [0,T]$ and all $x \in \mathbb{X}$. %
    Correspondingly, we define the inflow boundary (depending on~$u$) of the PDE via 
    $$
        \partial \XX^\text{in}(t) := \{x\in \partial \XX\,|\,a(t,x)\cdot\nu(x) > 0 \}.
    $$ 
    We now formulate two auxiliary variational problems to analyze the projection error. First, we consider for all $(w,v)\in L^2(\mathbb{X})\times L^2(\partial\XX)$ and $t\in (0,T)$,
    \begin{align}\label{eq:transport_boundary}
    \begin{split}
        \tfrac{\mathrm{d}}{\mathrm{d}t} \langle \Phi(t),w\rangle_{L^2(\mathbb{X)}} &= \langle a(t,\cdot)\cdot\nabla \Phi(t),w \rangle_{L^2(\mathbb{X})}\\
        \langle \Phi(t),v\rangle_{L^2(\partial \XX^\text{in}(t))} & = \langle \varphi,v\rangle_{L^2(\partial \XX^\text{in}(t))} \\
        \langle \Phi(0),w\rangle_{L^2(\mathbb{X})} &= \langle \varphi,w\rangle_{L^2(\mathbb{X})}.
    \end{split}
    \end{align}
    As the boundary values on the inflow boundary $\partial \XX^\text{in}(t)$ are prescribed, this transport equation is well-posed \cite[Chapter 14]{QuarVall2008}. Moreover, it can be straightforwardly verified that $\Phi\in \mathcal{C}(0,T;H^1(\XX))$ defined by
    \begin{align}\nonumber 
        \Phi(t)(x_0)=\begin{cases}
            \varphi(x(t;x_0,u)) &\text{if } x(s;x_0,u)\in \XX\,\,\forall s\in [0,t]\\
            \varphi(x_\text{exit}) &\text{otherwise}
        \end{cases}
    \end{align}
    solves \eqref{eq:transport_boundary}, where $x_\text{exit} \in \partial \XX^{\text{in}}(t)$ is the point at which $x(t;x_0,u)$ leaves~$\XX$.

Similarly, we consider the projected system such that for all test functions $(w_\mathbb{V},v_\mathbb{V})\in \mathbb{V}\times\mathbb{V}_\partial$, where $\mathbb{V}_\partial$ consists of the traces of the observable functions in $\mathbb{V}$, and $t\in (0,T)$,
\begin{align}\label{eq:transport_boundary_proj}
\begin{split}
        \tfrac{\mathrm{d}}{\mathrm{d}t} \langle \Phi_\mathbb{V}(t),w_\mathbb{V}\rangle_{L^2(\mathbb{X})} &= \langle a(t,\cdot)\cdot \nabla \Phi_\mathbb{V}(t),w_\mathbb{V}\rangle_{L^2(\mathbb{X})}\\
        \langle \Phi_\mathbb{V}(t),v_\mathbb{V}\rangle_{L^2(\partial \XX^\text{in}(t))} &= \langle \varphi,v_\mathbb{V}\rangle_{L^2(\partial \XX^\text{in}(t))} \\
    \langle \Phi_\mathbb{V}(0),w_\mathbb{V}\rangle_{L^2(\mathbb{X})} &= \langle \varphi,w_\mathbb{V}\rangle_{L^2(\mathbb{X})},
\end{split}
\end{align}
whose solution~$\Phi_\mathbb{V}$ is given by the projection of its counterpart for the variational problem~\eqref{eq:transport_boundary} 
onto $\mathbb{V}$.

The solutions of the variational problems~\eqref{eq:transport_boundary} and~\eqref{eq:transport_boundary_proj} coincide with the flow of the Koopman resp.\ the Koopman surrogate model on the set of initial values, such that the flow is contained in $\XX$. 
More precisely, for $z(t) = \varphi(x(t,\cdot;u))$ satisfying~\eqref{eq:bilin_real} and $z_\mathbb{V}(t)$ solving the surrogate dynamics~\eqref{eq:bilin_proj}, 
we have 
\begin{align}\label{eq:coincide}
    \Phi(t)(x_0) = z(t)(x_0) \quad\text{and}\quad \Phi_\mathbb{V}(t)(x_0) = z_\mathbb{V}(t)(x_0)
\end{align}
for all $x_0 \in \mathcal{A}(t)$ and $t\in [0,T]$. 
As $\varphi \in \mathcal{C}^2(\mathbb{X},\mathbb{R})$, the projection error between the auxiliary problems~\eqref{eq:transport_boundary} and~\eqref{eq:transport_boundary_proj}, i.e., the difference between $\Phi$ and $\Phi_\mathbb{V}$, can be bounded using finite element convergence results, cf.\ \cite[Section 14.3]{QuarVall2008}. In our case of linear finite elements, an application of \cite[Inequality~(14.3.16)]{QuarVall2008} reads
\begin{align}\label{eq:feestimate}
    \left( \int_{\XX} (\Phi(t)(x)- \Phi_\mathbb{V}(t)(x))^2\,\mathrm{d}x \right)^{1/2} \leq c \Delta x
\end{align}
for a constant $c=c(\|\varphi\|_{H^2(\XX)},|\mathbb{X}|,\underline{a},\overline{a})\geq 0$ and all $t\in [0,T]$.
Thus,
\begin{align*}
    \int_{\mathcal{A}(t)} ( \underbrace{ \varphi(x(t;\hat{x},u)) }_{ = z(t)(\hat{x}) } & - \tilde{z}_{m}(t)(\hat{x}) )^2\,\mathrm{d}\hat{x} \\
    &\leq  2 \int_{\mathcal{A}(t)}\! ( \underbrace{ z(t)(\hat{x}) - z_{\mathbb{V}}(t)(\hat{x})}_{ \stackrel{\eqref{eq:coincide}}{=} \Phi(t)(\hat{x})-\Phi_{\mathbb{V}}(t)(\hat{x}) } )^2 \!+\! (z_{\mathbb{V}}(t)(\hat{x}) - \tilde{z}_{m}(t)(\hat{x}) )^2\,\mathrm{d}\hat{x}.
\end{align*}
Taking square roots, the first term 
is bounded by $\nicefrac{\varepsilon}{2}$ for a mesh width $\Delta x=\mathcal{O}(\varepsilon)$ using~\eqref{eq:feestimate}.
The second term can be estimated by $\nicefrac{\varepsilon}{2}$ with probabilistic tolerance~$\delta$ using Corollary~\ref{cor:trajectories} with $\underbar{m}=\mathcal{O}\left(\nicefrac{N^2}{\varepsilon^2\delta}\right)$. 
Then, the result follows for dictionary size $N =\mathcal{O}\left(\nicefrac{1}{\Delta x^d} \right)=\mathcal{O}\left(\nicefrac{1}{\varepsilon^d} \right)$. 
\end{proof}

\begin{rem}
    On a $d$-dimensional domain~$\mathbb{X}$, Theorem~\ref{thm:fem} yields data requirements $\underbar{m} = \mathcal{O}(\varepsilon^{-2(d+1)})$ to approximate the generator and, thus, \textit{suffers} from 
    the curse of dimensionality, see also~\cite{ZhanZuaz21} for a comparison of eDMD for system identification to other methods. 
    Thus, augmenting the state by the control would exponentially scale the data requirements w.r.t.\ the input dimension, that is, $m=\mathcal{O}(\varepsilon^{-2(d+2+n_c)})$.
    In contrast, the proof of Corollary~\ref{cor:composite} reveals that the data requirements satisfy $m=\mathcal{O}((n_c+1)\varepsilon^{-(2(d+1))})$, i.e., linear scaling. 
\end{rem}

\section{Optimal and Model Predictive Control}\label{sec:control}
In this section, we show the usefulness 
of the derived uniform error bound in data-based optimal and predictive control. 

To this end, we consider the 
Optimal Control Problem 
\begin{align}\label{eq:OCP}
    \text{Minimize}&_{u \in \mathcal{U}_T(x_0)}\int_0^T %
    \ell(x(t;x_0,u),u(t))
    \,\mathrm{d}t \tag{OCP}
\end{align}
subject to 
the initial condition $x(0) = x_0$, the control-affine system dynamics~\eqref{eq:ode}, 
and the state constraints
\begin{align}\label{eq:constraints:state}
    h_j(x(t;x_0,u)) \leq 0 \qquad \forall\,j \in \{1,2,\ldots,p\}
\end{align}
for $t \in [0,T]$, %
where the set $\mathcal{U}_T(x_0)$ 
of admissible control functions is given by~\eqref{eq:admissible_controls}.
Further, we assume, that the set $\XX$ is chosen such that it contains the state constraint set in its interior, that is,
\begin{align*}
    \{x\in \mathbb{R}^{n} \,|\, h_j(x)\leq 0\,\,\text{for all }j \in \{1,2,\ldots,p\} \} \subsetneq \operatorname{int}(\mathbb{X}).
\end{align*}
The key challenge is to properly predict the performance index of~\eqref{eq:OCP} and ensure satisfaction of the state constraints~\eqref{eq:constraints:state} using 
the data-based surrogate model instead of propagating the state dynamics and then evaluating the \textit{observables} of interest, cp.~Figure~\ref{fig:sketch} and recall the identity
\begin{align}\label{eq:realkoop}
    (\mathcal{K}^t_u \varphi)(x_0) = \varphi(x(t;x_0,u)).
\end{align}
Since the Koopman operator~$\mathcal{K}^t_u$ is, in general, not known analytically, %
we resort to eDMD as outlined in Section~\ref{Sec:Recap} to derive a data-based finite-dimensional approximation~$\tilde{\mathcal{K}}^t_u$. %

All central quantities, i.e., the stage cost~$\ell$ and the constraint functions $h_j$, $j \in \{1,\ldots,p\}$,
are evaluated along the system dynamics~\eqref{eq:ode}. %
Hence, we use the observables $\varphi = h_j$, $j \in \{1,2,\ldots,p\}$, to ensure satisfaction of the state constraints. Assuming separability of
the stage cost
\begin{equation}\label{eq:separability}
    \ell(x,u) = \ell_1(x) + \ell_2(u),
\end{equation} 
we choose $\varphi = \ell_1$ as an observable while $\ell_2$ is at our disposal anyway. We point out that the assumed separability is typically the case. 
Otherwise, one can consider the coordinate functions as observables, i.e., $\varphi(x) = x_i$ for $i \in \{1,\ldots,d\}$, to evaluate~$\ell$. Theorem~\ref{thm:fem} allows to rigorously ensure constraint satisfaction and a bound $\varepsilon > 0$ on the approximation error w.r.t.\ the stage cost provided that the amount of data is sufficiently large and the finite-element dictionary is sufficiently rich. 
Consequently, the following result allows us to approximately solve the problem~\eqref{eq:OCP} using the derived eDMD-based, bilinear surrogate model with guaranteed constraint satisfaction and performance. 
\begin{prop}[State constraint and stage cost]\label{prop:ocp}
    Let Assumptions~\ref{as:main} and~\ref{ass:finite_elements} hold. %
    Further, suppose that $\ell_1, h_i\in \mathcal{C}^2(\mathbb{X},\mathbb{R})$, $i \in \{1,2,\ldots,p\}$. %
    Then, for error bound $\varepsilon > 0$, probabilistic tolerance $\delta \in (0,1)$, optimization horizon $T>0$, and all measurable control functions $u:[0,T]\to\mathbb{U}$ the following estimates hold:
    \begin{enumerate}
        \item Averaged probabilistic \textit{performance bound}, i.e., \\[-5mm]
    \end{enumerate}
    \begin{align*}
            \mathbb{P} \left(\|\ell(x(t;\cdot,u),u(t)) - \tilde{\ell}_{m}(t;\cdot,u)\|_{L^2(\mathcal{A}(t))} \leq \varepsilon \right) \geq 1-\delta.
    \end{align*}
    \begin{enumerate}
        \item [(2)] Averaged probabilistic \textit{state-constraint satisfaction} if the tightened constraint \\ $\frac{1}{\sqrt{|\mathcal{A}(t)}|} \int_{\mathcal{A}(t)} \tilde{h}_{i,m} (t;{x}_0,u)\,\mathrm{d}x_0 \leq -\varepsilon$ holds, i.e.,\\[-4mm]
    \end{enumerate}
    \begin{align*}
            \mathbb{P} \Big(\tfrac{1}{\sqrt{|\mathcal{A}(t)|}} \int_{\mathcal{A}(t)}h_i(x(t;\hat{x},u))\,\mathrm{d}\hat{x} \leq 0\Big) \geq 1 - \delta
    \end{align*}
    for $\tilde{\ell}_m(t;x_0,u) = \tilde{\ell}_{1,m}(t;x_0,u) + \ell_2(u(t))$ and all $ i\in \{1,\ldots,p\}$, where $\tilde{\ell}_{1,m}$, $\tilde{h}_{i,m}$, $i \in \{1,2,\ldots,p\}$, are predicted along the bilinear surrogate dynamics~\eqref{eq:bilin_surr} with $\tilde{\ell}_{1,m}(0;x_0,u)=P_{\mathbb{V}}\ell_1$ and $\tilde{h}_{i,m}(0;x_0,u)=P_{\mathbb{V}}h_i$, respectively, provided that the number of data points $m\geq \underbar{m}(\varepsilon,\delta)$ and the mesh size $\Delta x \leq \nicefrac{\varepsilon}{c}$, with $c=c(\|\ell_1\|_{H^2(\XX)},\|h\|_{H^2(\XX,\mathbb{R}^p)})$ are chosen in 
        according to Theorem~\ref{thm:fem}. In particular $m$ and $\Delta x$ can be determined independently of the chosen control~$u$.
\end{prop}
\begin{proof}
For the first assertion, i.e., the claim 
w.r.t.\ the stage cost, we invoke the assumed separability to compute
\begin{align*}
    \ell(x(t;x_0,u),u(t)) - \tilde{\ell}_{m}(t;x_0,u) = \ell_1(x(t;x_0,u)) -  \tilde{\ell}_{1,m}(t;x_0,u).
\end{align*}
Hence, the claim follows by setting $\varphi = \ell_1$ in Theorem~\ref{thm:fem}. 

Next, we show the second claim. To this end, we set $\varphi = h_i$, $i \in  \{1,\ldots,p\}$, in Theorem~\ref{thm:fem} and use the Cauchy-Schwarz inequality to get
\begin{align*}
    \tfrac 1 {\sqrt{|\mathcal{A}(t)|}} &\int_{\mathcal{A}(t)} h_i(x(t;\hat{x},u)) - \tilde{h}_{i,m}(t;\hat{x},u)\,\mathrm{d}\hat{x}\leq  \|h_i(x(t;\cdot,u)) - \tilde{h}_{i,m}(t;\cdot,u)\|_{L_2(\mathcal{A}(t))} \leq \varepsilon.
\end{align*}
Then, invoking the assumption completes the proof by
\begin{align*}
    \int_{\mathcal{A}(t)}  \tilde{h}_{i,m}(t;\hat{x},u) & + (h_i(x(t;\hat{x},u)) -\tilde{h}_{i,m}(t;\hat{x},u))\,\mathrm{d}\hat{x} \\ 
    \leq & \int_{\mathcal{A}(t)}  h_i(x(t;\hat{x},u)) - \tilde{h}_{i,m}(t;\hat{x},u)\,\mathrm{d}\hat{x} - \sqrt{|\mathcal{A}(t)|} \varepsilon
    \leq 0. 
\end{align*}
\end{proof}
The error bound of Proposition~\ref{prop:ocp} is given in an average sense due to the $L^2$-bound in the projection error estimate \eqref{eq:feestimate} of Theorem~\ref{thm:fem}. The projection error vanishes if the dictionary $\mathbb{V}$ is invariant under the Koopman semigroup or equivalently the generator, e.g., if it is spanned by eigenfunctions, cf.~\cite{KKB21}. In this case, Proposition~\ref{prop:ocp} can be straightforwardly refined to ensure a pointwise bound w.r.t.\ the initial value due to Corollary~\ref{cor:trajectories}.

In view of Proposition~\ref{prop:ocp} bounding the stage cost error and yielding chance constraint satisfaction, we briefly provide an outlook with respect to predictive control.

\textbf{Towards Model Predictive Control}: OCPs also play a predominant role in optimization-based control techniques like Model Predictive Control (MPC), where Problem~\eqref{eq:OCP} on an infinite-time horizon, i.e., $T = \infty$, is approximately solved by solving~\eqref{eq:OCP} at successive time instants~$i \delta$, $i \in \mathbb{N}_0$, on the prediction horizon $[i\delta, i\delta + T]$ subject to the current state 
as initial value,
see, e.g., the monographs~\cite{GrunPann17} 
and~\cite{CoroGrun20} w.r.t.\ MPC for continuous-time systems. Having obtained rigorous error estimates in view of optimal control, this paves the way of analyzing data-driven MPC schemes as proposed in~\cite{Peitz2020} and~\cite{KM18a} w.r.t.\ recursive feasibility or 
stability.

\section{Conclusion and outlook}\label{sec:conclusions}

Motivated by data-based surrogate modeling for optimal control problems with state constraints, we derived quantitative error estimates for eDMD-approximations of control systems. In this context, we provided a novel bound for the estimation uniform in the control and generalized the error analysis of the projection error to control systems. Further, using these probabilistic bounds, we derived error bounds on the performance and satisfaction of state constraints in data-based optimal and predicted control.

In future work, we further elaborate the presented results towards optimal control to derive suboptimality estimates~\cite{CoroGrun20} depending on both data and dictionary size. Moreover, a sensitivity analysis of the OCP could reveal robustness of optimal solutions w.r.t.\ approximation errors, that can be further exploited by numerical techniques, cf.\ \cite{GrunScha21}. Furthermore, a comparison to other approximation techniques for the Koopman operator, e.g., based on neural networks as proposed by~\cite{WangLou22}, might be of interest.

\bibliographystyle{abbrv}
\bibliography{references.bib}

\end{document}